\documentclass[11pt]{article}
\usepackage[T1]{fontenc}
\usepackage{lmodern}
\usepackage[margin=1in]{geometry}
\usepackage{amsmath,amsfonts,amssymb,amsthm,mathtools}
\usepackage{microtype,booktabs,enumitem}
\usepackage[colorlinks=true,linkcolor=blue,citecolor=blue,urlcolor=blue]{hyperref}
\hypersetup{pdftitle={Biased orientations of random graphs},pdfauthor={Bojan Mohar}}

\title{Spectral distribution for biased orientations of random graphs}

\author{Bojan Mohar \\ Simon Fraser University\\ Burnaby, BC, Canada \\ and \\ FMF, University of Ljubljana\\ Ljubljana, Slovenia \\ \texttt{mohar@sfu.ca} \thanks{Supported in part by NSERC Discovery Grant R832714 (Canada), by the ERC Synergy grant (European Union, ERC, KARST, project number 101071836), and by the Research Project N1-0218 of ARIS (Slovenia).}
}
\date{}

\newcommand \DD{\mathcal D}
\newcommand \TT{\mathcal T}
\newcommand \RR{\mathbb R}
\newcommand \CC{\mathbb C}
\newcommand \EE{\mathbb E}
\newcommand \PP{\mathbb P}
\newcommand \scm{\mu_{\mathrm{sc}}}
\newcommand \scp{\rho_{\mathrm{sc}}}
\newcommand \op{\mathrm{op}}
\newcommand \wprob{\xrightarrow{\mathrm{w},\,\PP}}
\newcommand \prob{\xrightarrow{\PP}}
\DeclareMathOperator{\Tr}{Tr}
\DeclareMathOperator{\rank}{rank}
\DeclareMathOperator{\Var}{Var}
\newtheorem{theorem}{Theorem}
\newtheorem{lemma}{Lemma}
\newtheorem{corollary}{Corollary}
\newtheorem*{bernstein}{Two-sided matrix Bernstein inequality}
\numberwithin{equation}{section}
\setlist[enumerate]{label=\textup{(\roman*)},leftmargin=*,itemsep=5pt}
\begin{document}

\maketitle

\begin{abstract}
We study the Hermitian adjacency matrix of the second kind for biased
orientations of Erd\H{o}s--R\'enyi random graphs. In the model
$\DD(n,p,p')$, an edge is present with probability $p$ and, conditional
on its presence, is directed from its smaller to its larger endpoint
with probability $p'$. Both parameters may depend on $n$, and biased
tournaments are the special case $\DD(n,1,p')$. We prove a global
semicircle law whenever
$n\bigl[p(1-p)+3p^2p'(1-p')\bigr]\to\infty$, for both the centered
and the uncentered matrix, and also for the expected empirical
spectral measure. Our second theorem gives the exact spectrum of the
mean matrix and the limits of the extreme eigenvalues on the scale
$np$, under a sufficient density assumption. The main difference 
from the unbiased random orientations (where $p'=\tfrac{1}{2}$)
is emergence of multiple large eigenvalues.
The proofs include the
closed-walk counting, control of the moment variances, an explicit
diagonalization of the deterministic bias, and the passage from
centered to uncentered spectral measures.
\end{abstract}

\section{Introduction}

Let $G(n,p)$ be the Erd\H{o}s--R\'enyi random graph on the ordered
vertex set $[n]=\{1,\ldots,n\}$. We obtain a random oriented graph
by independently directing every present edge from its smaller to
its larger endpoint with probability $p'$, and in the opposite
direction with probability $1-p'$. We denote the resulting model
by $\DD(n,p,p')$. The parameters have distinct roles: $p$ controls
edge presence, whereas $p'$ controls orientation bias. Throughout,
both $p=p(n)$ and $p'=p'(n)$ may vary with $n$.

Biased tournament models have some history. Probably the earlies paper is by
Frank \cite{Frank1968}. {\L}uczak, Ruci\'nski, and Gruszka \cite{LuczakEtAl1996}
studied the evolution of biased tournaments, and. more recently, Gishboliner, Krivelevich,
and Kronenberg \cite{GishbolinerEtAl2018} considered the chromatic number of
biased tournaments.
However, their eigenvalues were not studied before. 
The same holds for biased orientations of random graphs.

We use the sixth-root Hermitian adjacency matrix introduced by
the author~\cite{Mohar2020}, also called the Hermitian adjacency matrix
of the second kind. There are several reasons why this specific choice
is the most natural one, see \cite{Mohar2020} and \cite{AkbariEtAl2025}.  
An arc has weight
$\omega=(1+i\sqrt3)/2$ in its direction and weight
$\overline\omega$ in the reverse direction. The earlier Hermitian
adjacency matrix, introduced independently by Guo and
Mohar~\cite{GuoMohar2017} and Liu and Li~\cite{LiuLi2015}, uses
weights $i$ and $-i$. Further properties of the second-kind matrix
are studied by Li and Yu~\cite{LiYu2022}. The main distinction between
the two conventions is that the real part
$\textrm{Re}(\omega)=1/2$ creates a nonzero mean even when $p'=1/2$.

Our results use Wigner's semicircle law~\cite{Wigner1958} and the
standard spectral comparison methods described for example in Bai and
Silverstein~\cite{BaiSilverstein2010} and in other textbooks 
about eigenvalues of random matrices, a popular example being a recent 
manuscript of Tao \cite{Tao2012}. 
Related results for random
oriented graphs and random mixed graphs include the work of
Chen, Li and Lian~\cite{ChenLiLian2013}, Hu, Li, Liu and
Zhang~\cite{HuEtAl2017}, and Guan et al.~\cite{GuanEtAl2022}.
Separation of the leading eigenvalue for first-kind Hermitian
matrices of random mixed graphs is studied by Hu, Broersma, Hou
and Zhang~\cite{HuEtAl2022}. The same authors give bounds on skew
spectral radii of general random oriented graphs
in~\cite{HuEtAl2024}. These results provide context, but the mean
and variance computation for the present matrix and model bear 
notable differences.

We give two main theorems for the family $\DD(n,p,p')$.
The first describes the bulk spectrum and its expected empirical
measure. The second describes the deterministic mean spectrum and
the extreme eigenvalues on the scale $np$. Uniformly oriented
random graphs and biased tournaments are particular special
cases, obtained by setting $p'=1/2$ and $p=1$, respectively.

The bulk theorem is a model-specific application of the classical
moment method. Its main additional ingredient is the treatment of
the mean matrix. Under a nonzero bias, the imaginary part of the
mean has large rank and norm of order $np$. Nevertheless, after
bulk normalization it can be decomposed into a matrix of rank
$o(n)$ and a matrix of norm $o(1)$. We derive this decomposition
from an explicit cotangent formula for its eigenvalues. Details
are included to explain both the parameter assumptions and the
different roles of bulk and extreme eigenvalues.

\subsection*{Comments about the use of AI}

I have known the results collected in this paper since the original
work in \cite{GuoMohar2017} and \cite{Mohar2020}. I have popularized
it in several talks since I like the $\tfrac{1}{2}$ shift used to 
deal with the skew-symmetric part. This shift gives the eigenvalues
at the complex line $\mathrm{Im}(z)=-\tfrac{1}{2}$, which is 
reminiscent to the location of zeros of the Riemann zeta function.
I also like the appearance of multiple large eigenvalues in the biased orientation
models.

I was saving the final technicalities and subsequent paper submission
for a graduate student. However, 
since the recent AI work in mathematics, I have realized that this
is no longer a good Master Thesis problem and have decided to write it myself.
I have used help from ChatGPT (the current model GPT-6.1 Sol) which
has enhanced my original draft manuscript to my satisfaction. 
The model did not go much beyond my current understanding, but it was
an extremely useful tool that saved me a lot of technical work and
also provided some more generality. For example, 
Sections \ref{sect:different distribution}--\ref{sect:bulk and exceptional}
and several comments scattered around have been the work of ChatGPT. 

My original results started with the unbiased case, but since that case has
been more or less covered by Hu et al.~\cite{HuEtAl2017,HuEtAl2022}
and by Guan et al.~\cite{GuanEtAl2022}, 
I decided to stick with biased orientations, the case of which provides 
a novel outcome, where several large eigenvalues appear in addition to the 
bulk of eigenvalues following the semicircle distribution.

\section{The model, normalization, and spectral measures}

\subsection{The random graph and its Hermitian matrix}

All graphs considered here have no loops, and there is at most one
arc between any unordered pair of vertices. More explicitly, for
each $u<v$, independently of all other unordered pairs, the model
$\DD(n,p,p')$ makes the choices
\begin{center}
\begin{tabular}{|l|c|}
\toprule
Outcome on $uv$ & Probability\\
\midrule
No edge & $1-p$\\
The arc $u\to v$ & $pp'$\\
The arc $v\to u$ & $p(1-p')$\\
\bottomrule
\end{tabular}
\end{center}
This is exactly the two-stage construction of first sampling
a random undirected graph
$G(n,p)$ and then orienting its present edges with bias $p'$.
In particular,
\begin{equation}\label{eq:tournamentcase}
 \TT(n,p'):=\DD(n,1,p')
\end{equation}
is the biased tournament on the ordered
vertex set $[n]$, with forward orientation probability $p'$. At
$p'=0$ or $1$ every present edge is directed consistently with
one of the two vertex orders. Such an orientation is deterministic
\emph{conditional on the underlying graph}; the graph remains
random when $0<p<1$.

For an oriented graph $D$, its Hermitian adjacency matrix is
$H=H(D)\in\CC^{V(D)\times V(D)}$, with
\begin{equation}\label{eq:Hdef}
 H_{uv}=\begin{cases}
 \omega,&u\to v\text{ is an arc},\\
 \overline\omega,&v\to u\text{ is an arc},\\
 0,&u=v\text{ or }u,v\text{ are nonadjacent},
 \end{cases}
\end{equation}
where
$$
 \omega=\frac{1+i\sqrt3}{2}.
$$
Thus $H_{vu}=\overline{H_{uv}}$, so $H=H^*$ is Hermitian and
has $n$ real eigenvalues. We count multiplicities and order them as
\[
 \lambda_1(H)\ge\lambda_2(H)\ge\cdots\ge\lambda_n(H).
\]
If $A=A(D)$ is the adjacency matrix of the underlying undirected graph
and $S=S(D)$ is the real skew-adjacency matrix, whose $uv$-entry is
$1$ for $u\to v$, $-1$ for $v\to u$, and $0$ otherwise, then
\begin{equation}\label{eq:decomposition}
 H(D)=\frac12 A(D)+i\,\frac{\sqrt3}{2}S(D).
\end{equation}
The two summands use the same edge choices and are generally
dependent. Our proofs work with the entries of $H$ directly.
Reversing every arc replaces $H$ by $\overline H$, which has the
same real eigenvalues. Consequently the spectral distributions
for parameters $p'$ and $1-p'$ are identical.

As mentioned before, several results work for more general \emph{gain graphs},
where each edge $u\to v$ has an arbitrary complex weight $\zeta$ with $|\zeta|=1$
and the symmetric matrix entry has $\overline\zeta$. The gain $\zeta$ may 
even be different for different pairs $u,v$. In particular, Theorem \ref{thm:bulk} 
extends to gain graphs.  

\subsection{The mean and variance of matrix entries}

For the remainder of the article, write
\begin{equation}\label{eq:parameters}
 \beta=2p'-1,\quad a=\frac{\sqrt3}{2}\beta,\quad
 \alpha=\frac12+i\,a,\quad
 \tau=|\alpha|^2=\frac14+a^2=1-3p'(1-p').
\end{equation}
We will first determine the mean and variance of the entries of $H$, viewed as a random matrix 
corresponding to our model $\DD(n,p,p')$.
Since $\tfrac{1}{4}\le\tau\le1$, the following expressions are nonnegative.
For $u<v$, direct calculation gives
\begin{align}
 \EE H_{uv}&=pp'\omega+p(1-p')\overline\omega=p\alpha,
     \label{eq:entrymean}\\
 \EE|H_{uv}|^2&=p,\\
 \nu:=\Var H_{uv}&=\EE|H_{uv}-p\alpha|^2
     =p-p^2\tau
     =p(1-p)+3p^2p'(1-p').
     \label{eq:variance}
\end{align}
For a complex random variable, variance here means
$\Var Z=\EE|Z-\EE Z|^2$, not $\EE(Z-\EE Z)^2$.
Notice that $0\le\nu\le p$. It vanishes precisely when $p=0$,
or when $p=1$ and $p'\in\{0,1\}$.

Let $J_n$ be the all-ones matrix, $I_n$ the identity, and $S_n$
the deterministic skew-symmetric matrix with entries $1$ above
the diagonal, $-1$ below it, and $0$ on the diagonal. Then the \emph{mean matrix} $M_n$ and the corresponding \emph{centered matrix} $W_n$ are defined as
\begin{equation}\label{eq:mean}
 M_n:=\EE H_n=\frac p2(J_n-I_n)+i\,paS_n,\qquad
 W_n:=H_n-M_n.
\end{equation}
The entries of $W_n$ above the diagonal are independent,
centered, and have common variance $\nu$. They satisfy the uniform
bound
\begin{equation}\label{eq:bounded}
 |(W_n)_{uv}|\le1+p|\alpha|\le2.
\end{equation}
We also define the normalization factor 
\begin{equation}\label{eq:s_n}
 s_n:=\sqrt{n\nu}.
\end{equation}

Let us observe that we allow $p=p_n$ and $p'=p'_n$ to depend on $n$.
Therefore, also the constants in \eqref{eq:parameters} and \eqref{eq:variance}
depend on $n$ in general.
When necessary we write $\nu_n,a_n,p_n,p'_n$ to emphasize dependence
on $n$. A normalization by $s_n$ is used only when $\nu_n>0$.

\subsection{Empirical and expected spectral measures}

For a Hermitian matrix $B$ of order $n$, we define its \emph{spectral measure} $\mu_B$ 
and its \emph{cumulative distribution function} $F_B$ as
\begin{equation}\label{eq:ESD}
 \mu_B=\frac1n\sum_{j=1}^n\delta_{\lambda_j(B)},\qquad
 F_B(x)=\mu_B(( -\infty,x]).
\end{equation}
Here $\delta_t$ is unit mass at $t$. For a random matrix $B$, its 
\emph{expected empirical spectral measure} is the deterministic measure
$\overline\mu_B$ that satisfies
\begin{equation}\label{eq:expectedmeasure}
 \overline\mu_B(U)=\EE\mu_B(U)
  =\frac1n\EE\#\{j:\lambda_j(B)\in U\}
\end{equation}
for every Borel set $U$. It is generally different from
$\mu_{\EE B}$, the spectral measure of the mean matrix.

The standard \emph{semicircle distribution} $\scm$ has density
\begin{equation}\label{eq:sc}
 \scp(x)=\frac{1}{2\pi}\sqrt{4-x^2}\,
                       \boldsymbol{1}_{[-2,2]}(x).
\end{equation}
Its odd moments vanish and its $2m$-th moment is
$C_m=(m+1)^{-1}\binom{2m}{m}$. 

For convergence of (random) probability measures, we write
$\mu_{B_n}\wprob\mu$ if, for every bounded continuous real
function $f$, $\int f\,d\mu_{B_n}\prob\int f\,d\mu$. That is, 
for every $\varepsilon>0$,
\[
\mathbb P\!\left(
\left|\int f\,d\mu_n-\int f\,d\mu\right|>\varepsilon
\right) \longrightarrow 0.
\]
Thus each integral against a fixed test function converges in probability to its prescribed limit.
The notation $\xrightarrow{\mathrm w}$ denotes ordinary weak
convergence of deterministic probability measures. We refer the reader to \cite{Billingsley1999} for further details.

For a Hermitian matrix,
$\|B\|_{\op}=\max_j|\lambda_j(B)|$ is its \emph{operator norm} 
(also called the \emph{spectral radius} of $B$).
We write $Z_n=O_{\PP}(r_n)$ when $Z_n/r_n$ is bounded in
probability, and $Z_n=o_{\PP}(r_n)$ when $Z_n/r_n\prob0$.

\section{Our main theorems}

We will prove two theorems. The first one establishes the weak convergence of spectral measures.

\begin{theorem}[The semicircle bulk and its expectation]\label{thm:bulk}
Let $D_n$ have distribution $\DD(n,p_n,p'_n)$, and let
$H_n,W_n,\nu_n,s_n$ be as in \eqref{eq:variance}--\eqref{eq:s_n}.
If
\begin{equation}\label{eq:bulkcondition}
 n\nu_n=n\bigl[p_n(1-p_n)+3p_n^2p'_n(1-p'_n)\bigr]\longrightarrow\infty,
\end{equation}
then
\begin{equation}\label{eq:bulklimits}
 \mu_{W_n/s_n}\wprob\scm,\qquad
 \mu_{H_n/s_n}\wprob\scm,\qquad
 \overline\mu_{H_n/s_n}\xrightarrow{\mathrm w}\scm.
\end{equation}
In particular, for any fixed $x<y$, we have
\begin{equation}
 \frac1n\#\{j:x s_n<\lambda_j(H_n)\le y s_n\}
       \prob\int_x^y\scp(t)\,dt   \label{eq:countprob}\\
\end{equation}
and
\begin{equation} \frac1n\EE\#\{j:x s_n<\lambda_j(H_n)\le y s_n\}
       \longrightarrow\int_x^y\scp(t)\,dt.  \label{eq:countexpect}
\end{equation}
\end{theorem}

This is a global spectral statement. A proportion tending to $1$
of the eigenvalues lies in any fixed enlargement of
$[-2s_n,2s_n]$, but $o(n)$ eigenvalues can be distributed differently
and some individual eigenvalues can be much larger.
The next theorem describes such eigenvalues on the scale $np_n$.

\begin{theorem}[The mean spectrum and the extreme eigenvalues]
\label{thm:extreme}
Let $D_n$ have distribution $\DD(n,p_n,p'_n)$, and put
$c_n=|a_n|=(\sqrt3/2)|2p'_n-1|$.
\begin{enumerate}
\item If $c_n>0$, set $\theta_n=\arctan(2c_n)$. The eigenvalues
of $M_n=\EE H_n$, in decreasing order, are exactly
\begin{equation}\label{eq:meanspectrum}
 \lambda_j(M_n)=p_n\left[
 c_n\cot\left(\frac{\theta_n+(j-1)\pi}{n}\right)-\frac12\right],
 \qquad 1\le j\le n.
\end{equation}
If $c_n=0$, they are $(n-1)p_n/2$ once and $-p_n/2$ with
multiplicity $n-1$. If $p_n=0$, all eigenvalues are zero.
\item Let $\Delta_n=\sqrt{n\nu_n\log n}+\log n$. For all parameter sequences, we have
\begin{equation}\label{eq:uniformerror}
 \max_{1\le j\le n}|\lambda_j(H_n)-\lambda_j(M_n)|
       =O_{\PP}(\Delta_n).
\end{equation}
\item Suppose $np_n/\log n\to\infty$ and $p'_n\to r\in[0,1]$.
If $r\ne\tfrac{1}{2}$, define
$c_*=(\sqrt3/2)|2r-1|$ and $\theta_* =\arctan(2c_*)$.
For every fixed integer $k\ge1$,
\begin{align}
 \frac{\lambda_k(H_n)}{np_n}
       &\prob\frac{c_*}{\theta_*+(k-1)\pi},       \label{eq:positive}\\
 \frac{\lambda_{n+1-k}(H_n)}{np_n}
       &\prob-\frac{c_*}{k\pi-\theta_*}.          \label{eq:negative}
\end{align}
If $r=\tfrac{1}{2}$, then
\begin{equation}\label{eq:fair}
 \frac{\lambda_1(H_n)}{np_n}\prob\frac12 \quad \textrm{and~} \quad
 \frac{\max_{2\le j\le n}|\lambda_j(H_n)|}{np_n}\prob0.
\end{equation}
\end{enumerate}
\end{theorem}

It is worth noting that Theorem \ref{thm:extreme} does not require 
\eqref{eq:bulkcondition}, although its part (iii) uses a related 
assumption that $np_n/\log n\to\infty$.
The condition $np_n/\log n\to\infty$ is sufficient (but not necessary) for the stated
extreme-eigenvalue limits; no optimal threshold is asserted.
The scale $np_n$ permits $p_n\to0$. If instead $p_n\to p_*>0$,
the limits on the scale $n$ are obtained by multiplying the
right-hand sides of \eqref{eq:positive}--\eqref{eq:fair} by $p_*$.
For a fixed nonzero limiting bias, each fixed index from either
end has a nonzero limit on the scale $np_n$. This does not mean
that a positive proportion of the eigenvalues has that scale.

\section{Ingredients of the proofs}

\subsection{Exact diagonalization of the deterministic bias}

\begin{lemma}\label{lem:toeplitz}
For $b\ge0$ and $c>0$, let
$K_n(b,c)=b(J_n-I_n)+icS_n$ and
$\theta=\arg(b+ic)\in(0,\pi/2]$. Then
\begin{equation}\label{eq:toeplitz}
 \lambda_j(K_n(b,c))
   =c\cot\left(\frac{\theta+(j-1)\pi}{n}\right)-b,
 \qquad 1\le j\le n.
\end{equation}
In particular,
\begin{equation}\label{eq:Sspectrum}
 \lambda_j(iS_n)=\cot\left(\frac{(2j-1)\pi}{2n}\right),
 \qquad \|iS_n\|_{\op}=\cot\frac{\pi}{2n}.
\end{equation}
Replacing $c$ by $-c$ leaves the spectrum unchanged.
\end{lemma}

\begin{proof}
Add $bI_n$ and write $A=bJ_n+icS_n$. Any eigenvalue $\gamma$ of
$A$ is real, and its eigenvector $z=(z_1,\ldots,z_n)$ satisfies
\begin{equation}\label{eq:row}
 \gamma z_u=b\sum_{v=1}^n z_v
       +ic\left(\sum_{v>u}z_v-\sum_{v<u}z_v\right).
\end{equation}
In the difference between row $u+1$ and row $u$, the first sum
cancels and the second expression changes by $-ic(z_{u+1}+z_u$.
Thus
\[
 (\gamma+ic)z_{u+1}=(\gamma-ic)z_u.
\]
Since $c\ne0$, $\gamma+ic$ is nonzero. If $z_1=0$, this
recurrence would force every coordinate to vanish, which is
impossible since eigenvector is non-zero. We may therefore normalize $z_1=1$ and write
\begin{equation}\label{eq:ratio}
 z_u=r^{u-1},\qquad r=\frac{\gamma-ic}{\gamma+ic}.
\end{equation}
Here $|r|=1$ and $r\ne1$. The first row of \eqref{eq:row} becomes
\[
 \gamma+ic=(b+ic)\sum_{u=1}^n r^{u-1}
            =(b+ic)\frac{1-r^n}{1-r}.
\]
Since $1-r=2ic/(\gamma+ic)$, cancellation gives
\begin{equation}\label{eq:boundary}
 r^n=\frac{b-ic}{b+ic}=e^{-2i\theta}.
\end{equation}
The $n$ roots are $r_j=e^{-2it_j}$, where
$t_j=(\theta+(j-1)\pi)/n\in(0,\pi)$.
Solving \eqref{eq:ratio} for $\gamma$ gives
\[
 \gamma=i\,c\,\frac{1+r_j}{1-r_j}=c\cot t_j.
\]
Conversely, these values and their geometric-coordinate vectors
determined by \eqref{eq:ratio}
satisfy the first row and every consecutive row difference;
hence they satisfy all rows. The $n$ values are distinct and
decreasing because cotangent is strictly decreasing on $(0,\pi)$.
They give the full spectrum of $A$. Subtracting $b$ proves
\eqref{eq:toeplitz}.

Taking $b=0$ and $c=1$ gives $\theta=\pi/2$ and yields
\eqref{eq:Sspectrum}. Finally,
$K_n(b,-c)=\overline{K_n(b,c)}$. Conjugating an eigenvector
preserves its real eigenvalue. This proves the last assertion.
\end{proof}

\subsection{The centered semicircle law by closed walks}

The next lemma uses standard trace formula to express spectral moments
$\sum_{j=1}^{n}\lambda_j^k = \Tr(X^k)$ via counting of closed walks in
the complete graph with edges weighted by the matrix entries. 
The convenient assumption is that the random matrix $X$ has centered entries $X_{uv}$.

The lemma allows the entry distribution and its variance to
depend on $n$. The entries need not have a rotationally symmetric
complex distribution, and no condition on $\EE X_{uv}^2$ is needed.
In particular, they can arise by ``centering'' any gain values in $H$,
not necessarily restricted to the sixth roots of unity.

\begin{lemma}\label{lem:sc}
Let $X_n=(X_{uv})$ be Hermitian with zero diagonal and independent
entries above the diagonal. Suppose that for every $u<v$, we have
\[
 \EE X_{uv}=0,\quad \EE|X_{uv}|^2=\nu_n>0,
 \quad \textrm{and\/} \quad |X_{uv}|\le L,
\]
where $L$ is a fixed constant. Let $B_n=X_n/\sqrt{n\nu_n}$. If $n\nu_n\to\infty$, then
\[
 \mu_{B_n}\wprob\scm\quad \textrm{and~} \quad
 \overline\mu_{B_n}\xrightarrow{\mathrm w}\scm.
\]
\end{lemma}

\begin{proof}
We give the moment argument in five steps.

\smallskip\noindent\emph{Step 1: Expand the trace as a sum of walks.}
For a fixed positive integer $k$,
\begin{equation}\label{eq:moment}
 \EE\frac1n\Tr B_n^k
 =\frac1{n(n\nu_n)^{k/2}}
   \sum_{u_1,\ldots,u_k=1}^n
    \EE\left[X_{u_1u_2}\cdots X_{u_ku_1}\right].
\end{equation}
Each term corresponds to a closed walk
$u_1u_2\cdots u_ku_1$. Since the diagonal of $B_n$ is zero, we may
assume that any consecutive vertices $u_j, u_{j+1}$ are distinct. 
If an unordered edge $uv$ occurs just
once, its associated centered entry is independent of all other
factors, so this term also vanishes in the expectation.

For any remaining walk, let $e$ be its number of distinct unordered
edges and $t$ its number of distinct vertices. Every edge is used
at least twice, hence $e\le k/2$. The graph traced by the walk
is connected, so $t\le e+1$.

\smallskip\noindent\emph{Step 2: Bound on nonleading walk patterns.}
A walk pattern records the walk after relabeling vertices in
their order of first appearance. In other words, the pattern is obtained from $u_1u_2\dots u_k$
by omitting the repeated appearances of vertices and keeping only the first one.
There are only finitely many
such patterns for fixed $k$, and a pattern with $t$ vertices
has at most $n^t$ assignments of actual vertex labels.
If an edge occurs $m\ge2$ times, in either direction, its factor
in the expectation has absolute value at most
\[
 \EE|X_{uv}|^m\le L^{m-2}\EE|X_{uv}|^2=L^{m-2}\nu_n.
\]
Independence across unordered edges now bounds the expectation
for the walk by $L^{k-2e}\nu_n^e$. Thus a fixed pattern contributes
at most a constant depending on $k,L$ multiplied by
\begin{equation}\label{eq:walkbound}
 \frac{n^t \nu_n^e}{n(n\nu_n)^{k/2}}
        =n^{t-e-1}(n\nu_n)^{e-k/2}.
\end{equation}
If $e<k/2$, the second factor tends to zero and the first is at
most one. If $e=k/2$ but $t<e+1$, the first factor is at most
$n^{-1}$ and hence (\ref{eq:walkbound}) also tends to $0$. Therefore only the case
$k=2m$, $e=m$, $t=m+1$ can have a nonzero limit. The underlying graph of any such pattern
is a tree, and the walk traverses each of its edges exactly twice. Such patterns are
called \emph{tree patterns}.

\smallskip\noindent\emph{Step 3: Count the leading tree walks.}
For a tree pattern, each of its edges is
traversed exactly twice. Deleting any tree edge separates the
support into two parts. A closed walk crosses this cut equally
often in both directions, so the two traversals must be one in
each direction. Each edge consequently contributes
$\EE(X_{uv}X_{vu})=\EE|X_{uv}|^2=\nu_n$.
This is why the complex phases do not affect the limiting law.
Note also that the value in (\ref{eq:walkbound}) is equal to $1$ 
for every tree walk. Therefore, it suffices to count the number of tree patterns
and tree walks. 

For clarity, the number of tree walk patterns is the
Catalan number $C_m$. Mark the first traversal of each edge as
an up-step and the return traversal as a down-step. On a tree
with each edge used twice, a return is to the parent of the
current vertex, and an up-step discovers a new vertex. The
result is a path with $m$ up-steps and $m$ down-steps, never
below its starting height. Conversely, such a path reconstructs
a unique tree walk by introducing the next unused vertex
at every up-step and returning to the parent at every down-step.
The reflection argument counts these paths as
\[
 \binom{2m}{m}-\binom{2m}{m+1}
       =\frac1{m+1}\binom{2m}{m}=C_m.
\]
Each pattern has $n(n-1)\cdots(n-m)$ distinct label assignments.
Its normalization in \eqref{eq:moment} is $n^{m+1}\nu_n^m$,
and its entry expectation is $\nu_n^m$. We obtain
\begin{equation}\label{eq:momentlimit}
 \EE\frac1n\Tr B_n^k\longrightarrow
  \begin{cases}0,&k\text{ odd},\\ C_m,&k=2m.\end{cases}
\end{equation}

\smallskip\noindent\emph{Step 4: The spectral moments concentrate.}
We have determined the limit of expectation of $n^{-1}\Tr B_n^k$ 
($n\to\infty$) for every fixed $k$. Now, we need to analyse its
variance. 
Write $P_w$ for the product of entries associated with a walk $w$.
The variance of $n^{-1}\Tr B_n^k$ is a normalized sum of
$\EE(P_wP_{w'})-\EE P_w\EE P_{w'}$ over pairs $w,w'$ of closed walks
of length $k$. If their unordered edge-sets are disjoint, the
products are independent and the covariance is zero, even if
the walks share a vertex. Otherwise their union is connected.
If some union edge occurs just once among the $2k$ factors,
centering kills both terms of the covariance. We therefore
only need unions with $e\le k$ edges, each used at least twice,
and $t\le e+1$ vertices.

The joint expectation is bounded by $L^{2k-2e}\nu_n^e$ (see Step~2).
The product of the separate expectations satisfies the same
bound whenever it is nonzero: every shared edge supplies an
extra factor $\nu_n$, which can be bounded by $L^2$.
For each union pattern, the normalization is
$n^2(n\nu_n)^k$ and there are at most $n^t$ labelings. Its absolute
contribution is at most a constant times
\[
 n^{t-e-2}(n\nu_n)^{e-k}\le n^{-1},
\]
where we used that $n\nu_n\ge1$ for all sufficiently large $n$.
There are finitely many pairs of patterns for fixed $k$, so
\begin{equation}\label{eq:momentvariance}
 \Var\left(\frac1n\Tr B_n^k\right)=O_{k,L}(n^{-1}).
\end{equation}
Chebyshev's inequality and \eqref{eq:momentlimit} show that
each moment converges in probability to its semicircle value.

\smallskip\noindent\emph{Step 5: Pass from moments to measures.}
We give the details because weak convergence does not, by itself,
imply convergence of moments of unbounded functions. Write
\[
 m_k(B_n)=\int x^k\,d\mu_{B_n}(x)=\frac1n\Tr B_n^k,
 \qquad m_k=\int x^k\,d\scm(x).
\]
The preceding steps show that $m_k(B_n)\prob m_k$ for every
integer $k\ge0$; for $k=0$ both sides are identically one.

First take an arbitrary subsequence of matrix orders. By further
extraction, choose orders $n_j$ within the subsequence such that
\begin{equation}\label{eq:diagonalextraction}
 \PP\left\{\max_{0\le k\le j}|m_k(B_{n_j})-m_k|>1/j\right\}
       \le2^{-j}.
\end{equation}
This is possible because each maximum involves only finitely many
moments, all converging in probability. The probabilities in
\eqref{eq:diagonalextraction} are summable. The first
Borel--Cantelli lemma therefore gives an event of probability $1$
on which $m_k(B_{n_j})\to m_k$ simultaneously for every integer
$k\ge0$. Independence between different matrix orders is not
needed.

Fix an outcome in that event and set $\eta_j=\mu_{B_{n_j}}$.
We now have a deterministic sequence of probability measures.
Since $\int x^2\,d\eta_j\to m_2=1$, there is a finite constant
$C$ such that, for every $R>0$,
\begin{equation}\label{eq:momenttightness}
 \sup_j\eta_j(\{|x|>R\})
       \le\frac1{R^2}\sup_j\int x^2\,d\eta_j(x)
       \le\frac C{R^2}.
\end{equation}
Thus the measures are tight. Every subsequence consequently has
a further subsequence converging weakly to a probability measure
$\eta$ on $\RR$. 

We next show that any such limit $\eta$ has the required moments.
To simplify the notation, we may assume that the originally chosen sequence 
converges weakly to $\eta$.
For a fixed $k$, choose an even integer $2q>k$. Convergence of
the $2q$-th moments gives
$K_{2q}:=\sup_j\int|x|^{2q}\,d\eta_j(x)<\infty$, and hence
\begin{equation}\label{eq:momenttails}
 \sup_j\int_{|x|>R}|x|^k\,d\eta_j(x)
       \le\frac{K_{2q}}{R^{2q-k}}\longrightarrow0
       \quad\text{as }R\to\infty.
\end{equation}
The weak limit $\eta$ also satisfies $\int|x|^{2q}\,d\eta\le K_{2q}$.
To see this, we consider the function $f_M(x):=\min\{|x|^{2q},M\}$ that is bounded and continuous, so the weak convergence 
$\eta_j \xrightarrow{\mathrm w} \eta$ implies 
\[
\int f_M\,d\eta = \lim_{j\to\infty}\int f_M\,d\eta_j \le K_{2q}.
\]
Letting $M\to\infty$ and applying monotone convergence proves the stated bound.

Let $\chi_R$ be a continuous function equal to one on $[-R,R]$,
zero outside $[-2R,2R]$, and taking values in $[0,1]$.
Weak convergence applies to the bounded continuous function
$x^k\chi_R(x)$. The above tail estimates allow us first to pass to
the weak limit for fixed $R$ and then let $R\to\infty$.
It follows that $\int x^k\,d\eta=m_k$ for every $k$.

These moments identify $\eta$ uniquely. Indeed, for $R>2$ and
every positive integer $q$,
\[
 \eta(\{|x|>R\})\le\frac{m_{2q}}{R^{2q}}
       =\frac{C_q}{R^{2q}}\le\left(\frac4{R^2}\right)^q.
\]
Letting $q\to\infty$ shows that $\eta$ is supported on $[-2,2]$.
It has the same integrals of all polynomials as the semicircle distribution $\scm$.
Uniform polynomial approximation of continuous functions on
$[-2,2]$ then implies that $\eta=\scm$. Thus, every weakly convergent
subsequence of $(\eta_j)$ has the same limit, and tightness
gives $\eta_j\xrightarrow{\mathrm w}\scm$.

We have proved that every original subsequence admits a further
subsequence along which weak convergence holds almost surely.
This implies weak convergence in probability. More explicitly,
if convergence failed for some bounded continuous test function,
there would be a subsequence along which the probability of a
fixed positive deviation stayed bounded away from zero. Its
almost surely convergent further subsequence would contradict
that property.

Finally, for bounded continuous $f$, put
$Y_n=\int f\,d\mu_{B_n}$ and $y=\int f\,d\scm$.
Then $Y_n\prob y$ and $|Y_n|,|y|\le\|f\|_\infty$. For any
$\varepsilon>0$,
\[
 \EE|Y_n-y|\le\varepsilon+
       2\|f\|_\infty\PP\{|Y_n-y|>\varepsilon\}.
\]
Taking limits and then letting $\varepsilon\downarrow0$ proves
convergence of the expectations, which is the weak convergence
$\overline\mu_{B_n} \xrightarrow{\mathrm w} \scm$.
\end{proof}

\subsection{Spectral comparison and a norm estimate}

For two Hermitian matrices $A,B$ of order $n$, the rank inequality
and Weyl's inequality state that
\begin{align}
 \sup_x|F_A(x)-F_B(x)|&\le\frac{\rank(A-B)}n,
                                                     \label{eq:rank}\\
 \max_j|\lambda_j(A)-\lambda_j(B)|&\le\|A-B\|_{\op}.
                                                     \label{eq:weyl}
\end{align}
See Bai and Silverstein~\cite{BaiSilverstein2010}.
Suppose $A=B+R+E$, with $\rank R\le r$ and $\|E\|_{\op}\le\varepsilon$.
The norm estimate moves each ordered eigenvalue by at most
$\varepsilon$, and the rank estimate changes a distribution
function by at most $r/n$. Combining them gives
\begin{equation}\label{eq:ranknorm}
 F_B(x-\varepsilon)-\tfrac{r}{n}
       \le F_A(x)\le F_B(x+\varepsilon)+\tfrac{r}{n}.
\end{equation}
In particular, if $r/n\to0$, $\varepsilon\to0$, and $\mu_B$
converges to a measure with continuous distribution function,
then $\mu_A$ has the same weak limit. The assertion also holds
for convergence in probability, by applying the displayed bounds
at fixed points and using continuity of the limiting function.

We will use the concentration inequality in the general form used
below. It is the two-sided consequence of
Tropp~\cite[Theorem~6.1(ii), pp.~421--422]{Tropp2012}; the
one-sided estimate is also stated as Theorem~1.4 of that paper.

\begin{bernstein}
Let $Y_1,\ldots,Y_N$ be independent random Hermitian matrices of
order $d$, with $\EE Y_k=0$ and $\|Y_k\|_{\op}\le L$ almost surely,
where $L>0$ is fixed.
Put
\[
 Z=\sum_{k=1}^N Y_k,\qquad
 \sigma^2=\left\|\sum_{k=1}^N\EE(Y_k^2)\right\|_{\op}.
\]
For every $t>0$,
\begin{equation}\label{eq:generalbernstein}
 \PP\{\|Z\|_{\op}\ge t\}
       \le2d\exp\left(-\frac{t^2}{2(\sigma^2+Lt/3)}\right).
\end{equation}
If $\sigma^2=0$, then $Z=0$ almost surely.
\end{bernstein}

\begin{proof}
Tropp's Theorem~6.1(ii) in \cite{Tropp2012} bounds
$\PP\{\lambda_{\max}(Z)\ge t\}$ by the right-hand side of
\eqref{eq:generalbernstein} with prefactor $d$.
Apply it to both $(Y_k)$ and $(-Y_k)$. The norm bound gives
$\lambda_{\max}(\pm Y_k)\le L$, and the variance parameter
is unchanged. Since
$\|Z\|_{\op}=\max\{\lambda_{\max}(Z),\lambda_{\max}(-Z)\}$,
a union bound supplies the factor $2$. 

If $\sigma^2=0$, positivity of $\EE(Y_k^2)$ implies $Y_k=0$ almost surely
for each $k$.
\end{proof}

To estimate the norm of $W_n$, write it as a sum of rank-1 matrices:
\begin{equation}\label{eq:pairsum}
 W_n=\sum_{u<v}Y^{(uv)},\qquad
 Y^{(uv)}=(W_n)_{uv}\, e_ue_v^*
           +\overline{(W_n)_{uv}}\, e_ve_u^*.
\end{equation}
Here $e_u$ is the $u$-th coordinate vector. The summands are
independent centered Hermitian matrices. The only nonzero block
of a summand has eigenvalues $\pm|(W_n)_{uv}|$, so its norm
is at most $2$ by \eqref{eq:bounded}. Squaring that block gives
\[
 \EE\bigl[(Y^{(uv)})^2\bigr]
      =\nu_n(e_ue_u^*+e_ve_v^*),\qquad
 \sum_{u<v}\EE\bigl[(Y^{(uv)})^2\bigr]=(n-1)\nu_n I_n.
\]
Apply \eqref{eq:generalbernstein} with $d=n$, $L=2$, and
$\sigma^2=(n-1)\nu_n$. For $t>0$ and $n\ge2$, this gives
\begin{equation}\label{eq:bernstein}
 \PP\{\|W_n\|_{\op}\ge t\}
 \le2n\exp\left(-\frac{t^2}
                 {2((n-1)\nu_n+2t/3)}\right).
\end{equation}
For example, take
$t=\sqrt{6(n-1)\nu_n\log n}+4\log n$.
The exponent has absolute value at least $3\log n$:
in general, $t=\sqrt{2\sigma^2 u}+2Lu/3$ satisfies
$t^2\ge2u(\sigma^2+Lt/3)$, and here
$\sigma^2=(n-1)\nu_n$, $L=2$, $u=3\log n$.
Thus the probability is at most $2n^{-2}$. In particular,
\begin{equation}\label{eq:normbound}
 \|W_n\|_{\op}
   =O_{\PP}\bigl(\sqrt{n\nu_n\log n}+\log n\bigr).
\end{equation}
When $\nu_n=0$, $W_n=0$ almost surely, and the same conclusion
holds. This estimate is intended for separation on the scale
$np_n$; it does not locate the sharp bulk edge.

\section{Proofs of the main theorems}

\subsection{Proof of Theorem \ref{thm:bulk}}

\begin{proof}
By the assumption \eqref{eq:bulkcondition}, we have that $s_n=\sqrt{n\nu_n}\to\infty$.
The entries of $W_n$ satisfy the assumptions of Lemma~\ref{lem:sc} with $L=2$.
Consequently $\mu_{W_n/s_n}\wprob\scm$.

To prove the convergence of spectral measures of $H_n/s_n$,
we will need to control the spectrum of the imaginary mean
$Q_n=i\,d_nS_n$, where $d_n=p_na_n$. 
We will use Lemma \ref{lem:toeplitz} with $b=0$ and $c=d_n$.
We may assume that $p'>\tfrac{1}{2}$, so we have $c>0$.
By \eqref{eq:Sspectrum}, the multiset of 
eigenvalues of $Q_n$ consists of values $d_n\cot\vartheta_j$, where
$\vartheta_j=(j-\tfrac{1}{2})\pi/n$. For $t>0$ and $d_n\ne0$,
$|d_n\cot\vartheta_j|>t$ if and only if
$\vartheta_j<\phi$ or $\vartheta_j>\pi-\phi$, where $\phi=\arctan(|d_n|/t)$.
The angles $\vartheta_j$ are equally spaced by $\pi/n$. An interval of length $\phi$
contains at most $n\phi/\pi+1$ of them, so
\begin{equation}\label{eq:tailcount}
 \#\{j:|\lambda_j(Q_n)|>t\}
    \le\frac{2n}{\pi}\arctan\frac{|d_n|}{t}+2
    \le\frac{2n|d_n|}{\pi t}+2.
\end{equation}
If $d_n=0$, the left-hand side is zero and the bound is still valid.

Next, we split the mean into two ``small'' perturbations.
Take the spectral decomposition of $Q_n$ and retain in $R_n$
only eigenvalues of absolute value greater than $\sqrt{s_n}$
and the corresponding projections into their eigenspaces.
Let $E_n=Q_n-R_n$ contain the part of the spectral decomposition 
corresponding to the remaining eigenvalues. Then
$Q_n=R_n+E_n$, and by \eqref{eq:tailcount},
\[
 \frac{\rank R_n}{n}\le
       \frac{2|d_n|}{\pi\sqrt{s_n}}+\frac2n\to0,
 \quad \mathrm{and\/} \quad \|E_n\|_{\op}\le\sqrt{s_n},
\]
because $|d_n|\le\sqrt3/2$ and $s_n\to\infty$. Using the full mean in
\eqref{eq:mean}, we have the exact decomposition
\begin{equation}\label{eq:bulkperturb}
 \frac{H_n}{s_n}=\frac{W_n}{s_n}
       +\frac{R_n+(p_n/2)J_n}{s_n}
       +\frac{E_n-(p_n/2)I_n}{s_n}.
\end{equation}
The middle term has rank at most $\rank R_n+1=o(n)$.
The norm of the last term is at most
$s_n^{-1/2}+p_n/(2s_n)\to0$. The comparison
\eqref{eq:ranknorm} and continuity of the semicircle distribution
function now give $\mu_{H_n/s_n}\wprob\scm$.
This argument is needed because the norm of the entire mean
divided by $s_n$ need not tend to zero.

Finally, we take expectations and interval counts.
For every bounded continuous test function, convergence of its
integral in probability implies convergence of the expectation
by the boundedness argument at the end of Lemma~\ref{lem:sc}.
This proves the last assertion of \eqref{eq:bulklimits}.
The endpoints of $(x,y]$ have zero semicircle mass. Weak
convergence therefore gives \eqref{eq:countprob}. Its left-hand
side lies in $[0,1]$, so its expectation also converges, proving
\eqref{eq:countexpect}.
\end{proof}

\subsection{Proof of Theorem \ref{thm:extreme}}

\begin{proof}
(i) First, we identify the exact mean spectrum.
For $p_n>0$, \eqref{eq:mean} gives
\[
 M_n=p_n\left[\frac12(J_n-I_n)+i\,a_nS_n\right].
\]
If $a_n<0$, complex conjugation allows us to replace it by
$c_n=|a_n|$ without changing the eigenvalues. Apply
Lemma~\ref{lem:toeplitz} with $b=1/2$ and $c=c_n$.
Its angle is $\arg(1/2+ic_n)=\arctan(2c_n)$, which proves
\eqref{eq:meanspectrum}. If $a_n=0$, the all-ones vector is an
eigenvector of $J_n-I_n$ with eigenvalue $n-1$, and its orthogonal
complement has eigenvalue $-1$. This gives the stated spectrum.
The case $p_n=0$ is immediate.

(ii) For part (ii), we pass from the mean to the random matrix.
Weyl's inequality applied to $H_n=M_n+W_n$, followed by
\eqref{eq:normbound}, proves \eqref{eq:uniformerror}.

(iii) If $np_n/\log n\to\infty$, then $p_n>0$ eventually, and
$\nu_n\le p_n$ gives
\begin{equation}\label{eq:relativeerror}
 \frac{\Delta_n}{np_n}
 \le\sqrt{\frac{\log n}{np_n}}+\frac{\log n}{np_n}
 \longrightarrow0.
\end{equation}
Consequently, the eigenvalues of $H_n/(np_n)$ and $M_n/(np_n)$
differ uniformly by $o_{\PP}(1)$. It remains to compute
deterministic limits.

We can now evaluate the extreme deterministic eigenvalues.
Suppose first that $p'_n\to r\ne \tfrac{1}{2}$. Then
$c_n\to c_*>0$, $\theta_n\to\theta_*>0$, and for fixed $k$,
\[
 \frac{\lambda_k(M_n)}{np_n}
   =\frac{c_n}{n}\cot\frac{\theta_n+(k-1)\pi}{n}
        -\frac1{2n}
   \longrightarrow\frac{c_*}{\theta_*+(k-1)\pi},
\]
since $x\cot x\to1$ as $x\to0$. At the other end, the angle
for $j=n+1-k$ is
$\pi-(k\pi-\theta_n)/n$. The identity
$\cot(\pi-x)=-\cot x$ gives
\[
 \frac{\lambda_{n+1-k}(M_n)}{np_n}
   =-\frac{c_n}{n}\cot\frac{k\pi-\theta_n}{n}
        -\frac1{2n}
   \longrightarrow-\frac{c_*}{k\pi-\theta_*}.
\]
Together with the above, this gives \eqref{eq:positive} and
\eqref{eq:negative}.

If instead $p'_n\to1/2$, then $a_n\to0$. The norm formula
in \eqref{eq:Sspectrum} and $\cot x\le1/x$ for $0<x\le\pi/2$
give
\[
 \left\|\frac{M_n}{np_n}-\frac{J_n-I_n}{2n}\right\|_{\op}
   =\frac{|a_n|}{n}\cot\frac{\pi}{2n}
   \le\frac{2|a_n|}{\pi}\longrightarrow0.
\]
The comparison matrix has eigenvalues $(n-1)/(2n)$ once and
$-1/(2n)$ otherwise. Further application of Weyl's inequality
and the above comparison of eigenvalues of $H_n$ and $M_n$
proves \eqref{eq:fair}.
\end{proof}

\section{Special cases and interpretation}

\subsection{How the familiar families fit into the model}

The following table records the bulk variance and the sufficient
bulk condition used in the proof of Theorem~\ref{thm:bulk}.
\begin{center}
\begin{tabular}{|l|c|l|}
\toprule
Special case & Entry variance $\nu$ & Bulk condition\\
\midrule
$p'=1/2$ & $p-p^2/4$ & $np\to\infty$\\[2pt]
$p=1$ & $3p'(1-p')$ & $np'(1-p')\to\infty$\\[2pt]
$p'\in\{0,1\}$ & $p(1-p)$ & $np(1-p)\to\infty$\\
\bottomrule
\end{tabular}
\end{center}
Both parameters in this table may still vary with $n$, subject
to the indicated restriction. In the first row, the equivalence
of the displayed condition with $n\nu\to\infty$ follows from
$3p/4\le p-p^2/4\le p$.

For unbiased orientations, Theorem~\ref{thm:extreme} gives the
more explicit estimates
\begin{align}
 \lambda_1(H_n)&=\frac{(n-1)p_n}{2}+O_{\PP}(\Delta_n),\\
 \max_{j\ge2}|\lambda_j(H_n)+p_n/2|&=O_{\PP}(\Delta_n),
\end{align}
since the mean has only one eigenvalue outside its scalar part.
Thus $np_n/\log n\to\infty$ ensures a single leading
eigenvalue asymptotic to $np_n/2$, with every other eigenvalue
$o_{\PP}(np_n)$.

For tournaments, set $p_n=1$ throughout. The model is
$\TT(n,p'_n)=\DD(n,1,p'_n)$, the scale in the bulk theorem is
$\sqrt{3np'_n(1-p'_n)}$, and the scale in the extreme theorem
is $n$. Its density assumption holds automatically.
At $p'=1/2$, this agrees with the unbiased orientation of
$G(n,1)$, and both descriptions give variance $3/4$ and a
leading eigenvalue asymptotic to $n/2$.

If $p=1$ and $p'\in\{0,1\}$, there is no randomness. The
tournament is transitive, up to reversing all arcs. Now
$c=\sqrt3/2$ and $\theta=\pi/3$, so its entire spectrum is
\begin{equation}\label{eq:transitive}
 \lambda_j(H_n)=\frac{\sqrt3}{2}
   \cot\left(\frac{\pi/3+(j-1)\pi}{n}\right)-\frac12.
\end{equation}
In particular, $\lambda_1(H_n)/n\to3\sqrt3/(2\pi)$ and
$\lambda_n(H_n)/n\to-3\sqrt3/(4\pi)$.
This recovers the computation in \cite{AkbariEtAl2026}.
The extreme-eigenvalue formulas still apply, while the bulk
theorem has no variance normalization in this deterministic case.

More generally, if $p_n\le1-\delta$ for a fixed $\delta>0$,
then $\delta p_n\le \nu_n\le p_n$, regardless of $p'_n$.
The bulk condition is then equivalent to $np_n\to\infty$.
If instead $p'_n\in[\delta,1-\delta]$ with
$0<\delta\le\tfrac{1}{2}$, then
$\nu_n\ge3\delta(1-\delta)p_n$, and the same equivalence holds.

\subsection{Expected counts and the scale of the bulk}

Theorem~\ref{thm:bulk} answers the question of the expected
eigenvalue distribution through the expected empirical measure.
For fixed $x<y$, it gives
\begin{equation}\label{eq:expectedcounts}
 \EE\#\{j:x s_n<\lambda_j(H_n)\le y s_n\}
       =n\int_x^y\scp(t)\,dt+o(n).
\end{equation}
On the original eigenvalue scale, a convenient bulk density
approximation is therefore
\begin{equation}\label{eq:density}
 f_n^{\mathrm{bulk}}(x)
   =\frac{1}{2\pi s_n^2}\sqrt{(4s_n^2-x^2)_+}
\end{equation}
where $y_+=\max\{y,0\}$.
Its meaning is the global interval-count statement above.
Since $n$ is finite, the actual expected spectral measure is discrete;
no exact density or local eigenvalue-spacing law is asserted.
The scalar shift $-p_n/2$ in the mean is negligible relative
to $s_n\to\infty$.

Under the extreme theorem's density assumption, we have $np_n\to\infty$
and $s_n\le\sqrt{np_n}$. Whenever the bulk theorem also applies,
\[
 \frac{np_n}{s_n}\ge\sqrt{np_n}\longrightarrow\infty.
\]
Thus eigenvalues with a nonzero limit on the scale $np_n$ are
well beyond the bulk scale. Their number need not be bounded,
but their proportion tends to zero by the bulk law.

\subsection{The spectrum of the mean is a different distribution}
\label{sect:different distribution}

For fixed $p>0$ and fixed $p'\ne\tfrac{1}{2}$, let $c=|a|>0$.
The angles in \eqref{eq:meanspectrum} become uniformly
distributed on $(0,\pi)$. It follows that the unscaled mean
matrix has the weak limit
\begin{equation}\label{eq:Cauchy}
 \mu_{M_n}\xrightarrow{\mathrm w}\mathsf C_{p,c},\quad \textrm{where~} \quad
 \frac{d\mathsf C_{p,c}}{dx}(x)
    =\frac{pc}{\pi((x+p/2)^2+p^2c^2)}.
\end{equation}
Indeed, this is the distribution of $pc\cot U-p/2$, with
$U$ uniform on $(0,\pi)$. For any bounded continuous test
function, restrict first to $U\in[\varepsilon,\pi-\varepsilon]$,
where ordinary Riemann sums apply; the two omitted end intervals
have total probability $2\varepsilon/\pi$. Letting
$\varepsilon\downarrow0$ proves the weak convergence.
Changing variables $x=pc\cot u-p/2$ gives the displayed density.
It is a Cauchy law with location $-p/2$ and scale $pc$.

For $p'=1/2$, the mean instead has limiting measure
$\delta_{-p/2}$, since its exceptional leading eigenvalue has
mass only $1/n$. In the deterministic transitive case
$p=1$, $p'\in\{0,1\}$, \eqref{eq:Cauchy} is also the
limiting unscaled spectral measure of $H_n$ itself. For fixed
parameters with $\nu>0$, random fluctuations create the growing
bulk scale $\sqrt{n\nu}$. Consequently the measure in
\eqref{eq:Cauchy} must be distinguished from
$\overline\mu_{H_n}$.

\subsection{Why uncentered second moments need not be semicircular}

Weak convergence controls bounded continuous test functions,
and the function $x^2$ is unbounded. Large eigenvalues can
therefore contribute to second moments without changing the
weak bulk limit.
For every oriented graph with $m$ arcs,
\[
 \Tr H^2=\sum_{u,v}H_{uv}H_{vu}
          =\sum_{u,v}|H_{uv}|^2=2m.
\]
In our model, $\EE m=\binom n2p_n$, so, when $\nu_n>0$,
\begin{equation}\label{eq:secondmoment}
 \EE\int x^2\,d\mu_{H_n/s_n}(x)
    =\frac{n-1}{n}\cdot\frac{p_n}{\nu_n}.
\end{equation}
For fixed $p>0$ and fixed $p'$ with $\nu>0$, the limit is
\begin{equation}\label{eq:excess}
 \frac p\nu=\frac1{1-p\tau}=1+\frac{p^2\tau}{\nu}>1,
\end{equation}
whereas the semicircle second moment is $1$.
The extra term comes from the mean. To see this, observe that
$\Tr M_n^2=n(n-1)p_n^2\tau_n$, and
\[
 \EE\Tr H_n^2=\EE\Tr W_n^2+\Tr M_n^2,
\]
because the cross term has expectation zero. Dividing by
$n^2\nu_n$ explains both the semicircle contribution and the
excess. This is why the moment proof is applied to $W_n$,
followed by a separate comparison for $H_n$.

\subsection{Further questions}

The theorems describe global spectral measures and coarse
extreme-eigenvalue limits. Natural next questions are the
limiting laws when $n\nu_n$ stays bounded, sharp separation
thresholds for small $p_n$, and fluctuations of the extreme
eigenvalues. When $p'_n\to1/2$, the smaller bias eigenvalues
may interact with the bulk on intermediate scales. Comparing
the norm of $i\,p_na_nS_n$ with $s_n$ suggests examining
$|p'_n-1/2|$ of order $s_n/(np_n)$; for fixed $p>0$
this is the usual order of magnitude $n^{-1/2}$. The present results
do not determine that transition or its local statistics.

\subsection{The bulk and the large exceptional eigenvalues}
\label{sect:bulk and exceptional}

We finish by stating together what the two main theorems say
about the different parts of the spectrum. Here the \emph{global
bulk} means the eigenvalues described by the limiting empirical
measure on the scale $s_n=\sqrt{n\nu_n}$. For a fixed
$\varepsilon>0$, we use the window
$[-(2+\varepsilon)s_n,(2+\varepsilon)s_n]$ to make the assertion
that this part contains asymptotically almost all eigenvalues
precise. This window includes the spectral edges. In the
narrower terminology used for local spectral statistics, bulk
locations are normalized points in compact subintervals of
$(-2,2)$; neighborhoods of $-2$ and $2$ are the \emph{edges}.

\begin{corollary}[Semicircle bulk with an unbounded exceptional part]
\label{cor:bulkexceptional}
Suppose that
\[
 n\nu_n\to\infty,\quad \frac{np_n}{\log n}\to\infty, \quad \mathrm{and~}
 \quad p'_n\to r\in[0,1]\setminus\{\tfrac{1}{2}\}.
\]
Let $s_n=\sqrt{n\nu_n}$,
$c_*=(\sqrt3/2)|2r-1|$, and $\theta_* =\arctan(2c_*)$.
\begin{enumerate}
\item For every fixed $\varepsilon>0$,
\begin{equation}\label{eq:bulkpopulation}
 \#\{j:|\lambda_j(H_n)|\le(2+\varepsilon)s_n\}
       =n-o_{\PP}(n).
\end{equation}
Within this scale, the eigenvalue counts follow the semicircle
law: for every $-2<x<y<2$,
\[
 \frac1n\#\{j:x s_n<\lambda_j(H_n)\le y s_n\}
       \prob\int_x^y\scp(t)\,dt.
\]
In particular, every such interval contains a number of
eigenvalues proportional to $n$ with probability tending to one.
\item For every fixed positive integer $k$,
\begin{equation}\label{eq:exceptionalconstants}
 \frac{\lambda_k(H_n)}{np_n}\prob \Lambda_k,
 \qquad \frac{\lambda_{n+1-k}(H_n)}{np_n}\prob-\Lambda'_k,
\end{equation}
where
\[
 \Lambda_k=\frac{c_*}{\theta_*+(k-1)\pi}>0,\qquad
 \Lambda'_k=\frac{c_*}{k\pi-\theta_*}>0.
\]
Thus any fixed number of eigenvalues from each end has scale
$np_n$, which is asymptotically larger than the bulk scale.
\item For fixed $\varepsilon>0$, define
\begin{align*}
 N_n^+(\varepsilon)&=\#\{j:\lambda_j(H_n)>(2+\varepsilon)s_n\},\\
 N_n^-(\varepsilon)&=\#\{j:\lambda_j(H_n)<-(2+\varepsilon)s_n\}.
\end{align*}
For every fixed integer $K\ge1$,
\begin{equation}\label{eq:exceptionaldivergence}
 \PP\{N_n^+(\varepsilon)\ge K,\ N_n^-(\varepsilon)\ge K\}
       \longrightarrow1,
 \quad \mathrm{and~} \quad
 \frac{N_n^+(\varepsilon)+N_n^-(\varepsilon)}n\prob0.
\end{equation}
The exceptional counts therefore tend to infinity in probability
on both sides, while their total proportion tends to zero.
\end{enumerate}
\end{corollary}

\begin{proof}
Theorem~\ref{thm:bulk} gives the interval-count limit and
\eqref{eq:bulkpopulation}, since the semicircle law assigns mass
one to $[-2,2]$ and has no atoms. The constants in part (ii)
are those of Theorem~\ref{thm:extreme}(iii). Moreover,
$\nu_n\le p_n$ implies
\[
 \frac{np_n}{s_n}\ge\sqrt{np_n}\longrightarrow\infty.
\]
For each fixed $K$, part (ii) consequently yields
$\lambda_K(H_n)/s_n\prob+\infty$ and
$\lambda_{n+1-K}(H_n)/s_n\prob-\infty$.
By the ordering of the eigenvalues, with probability tending
to $1$, at least $K$ eigenvalues lie above the positive threshold and at
least $K$ lie below the negative threshold. This proves the
first assertion of \eqref{eq:exceptionaldivergence}. The second one
is precisely the complement of \eqref{eq:bulkpopulation}.
\end{proof}

The corollary distinguishes two types of conclusions. The bulk
is a statement about proportions and the empirical distribution.
The explicit exceptional limits concern fixed indices from
either end. They do not assert that every eigenvalue outside
the bulk window has order $np_n$. Indices growing with $n$ can
occur on intermediate scales. Nor do they specify a canonical
boundary between the bulk and its edges.

For comparison, consider the ordinary adjacency matrix $A_n$ of
$G(n,p)$ with fixed $p\in(0,1)$. Its mean is $p(J_n-I_n)$,
which, after a scalar shift, has rank one. The classical result
of F\"uredi and Koml\'os~\cite[Theorem~1]{FurediKomlos1981} gives
$\lambda_1(A_n)/(np)\prob1$ and, for every fixed
$\varepsilon>0$,
\begin{equation}\label{eq:ERsingleoutlier}
 \PP\left\{\max_{2\le j\le n}|\lambda_j(A_n)|
       \le(2+\varepsilon)\sqrt{np(1-p)}\right\}\longrightarrow1.
\end{equation}
Thus, with probability tending to one, there is exactly one
eigenvalue beyond any fixed enlargement of the ordinary
semicircle range. Under a nonzero limiting orientation bias,
Corollary~\ref{cor:bulkexceptional} instead gives unbounded
exceptional counts on both sides. The additional deterministic
term $ip_na_nS_n$ is responsible for this behavior: it has a
whole sequence of positive and negative eigenvalues of order
$np_n$ at fixed indices. At $p'_n=1/2$ this imaginary mean
vanishes, and Theorem~\ref{thm:extreme} recovers single exceptional
eigenvalue on the scale $np_n$. 

\bibliographystyle{plain}
\bibliography{biased_orientations}

\end{document}